\documentclass{amsart}

\usepackage{epsfig}

\usepackage{amsthm}

\usepackage{amssymb}

\usepackage{amsmath}

\usepackage{amscd}

\newtheorem {Theorem}  {Theorem}

\numberwithin{Theorem}{section}

\numberwithin{equation}{section}

\newtheorem {Lemma}[Theorem]  {Lemma}

\newtheorem {Proposition}[Theorem]{Proposition}

\theoremstyle{definition}

\newtheorem{Definition}[Theorem]{Definition}

\theoremstyle{remark}

\newtheorem{Remark}[Theorem]{Remark}

\newtheorem{Example}[Theorem]{Example}

\newtheorem {Corollary}[Theorem]{Corollary}

\expandafter\chardef\csname pre amssym.def at\endcsname=\the\catcode`\@ \catcode`\@=11

\def\undefine#1{\let#1\undefined}

\def\newsymbol#1#2#3#4#5{\let\next@\relax

 \ifnum#2=\@ne\let\next@\msafam@\else

 \ifnum#2=\tw@\let\next@\msbfam@\fi\fi

 \mathchardef#1="#3\next@#4#5}

\def\mathhexbox@#1#2#3{\relax

 \ifmmode\mathpalette{}{\m@th\mathchar"#1#2#3}%

 \else\leavevmode\hbox{$\m@th\mathchar"#1#2#3$}\fi}

\def\hexnumber@#1{\ifcase#1 0\or 1\or 2\or 3\or 4\or 5\or 6\or 7\or 8\or

 9\or A\or B\or C\or D\or E\or F\fi}

\font\teneufm=eufm10 \font\seveneufm=eufm7 \font\fiveeufm=eufm5

\newfam\eufmfam

\textfont\eufmfam=\teneufm \scriptfont\eufmfam=\seveneufm

\scriptscriptfont\eufmfam=\fiveeufm

\catcode`\@=\csname pre amssym.def at\endcsname

\newcounter{remark}

\newcommand{\rn}{{\mathbb R}^n}

\newcommand{\bg}{\begin{equation}}

\newcommand{\ed}{\end{equation}}

\newcommand{\bga}{\begin{eqnarray}}

\newcommand{\eda}{\end{eqnarray}}

\newcommand{\e}{\epsilon}

\newcommand{\gm}{\gamma}

\newcommand{\R}{\mathbf{R}}

\newcommand{\Aa}{{\mathcal A}}

\newcommand{\Bb}{{\mathcal B}}

\newcommand{\Vv}{{\mathcal V}}

\newcommand{\Ss}{{\mathcal S}}

\def  \C   {{\mathbb C}}

\def  \R   {{\mathbb R}}

\def  \12  {{\frac{1}{2}}}

\newcommand{\LL}{{\rm L}}

\newcommand{\HH}{{\rm H}}

\newcommand{\WW}{{\rm W}}

\begin{document}







\title[Poincar\'e's inequality and diffusive evolution equations]{Poincar\'e's inequality and diffusive evolution equations}

\author[Clayton Bjorland ]{Clayton Bjorland }

\address{Department of Mathematics, UC Santa Cruz, Santa Cruz, CA 95064, USA}


\author[Maria E. Schonbek ]{Maria E. Schonbek }

\address{Department of Mathematics, UC Santa Cruz, Santa Cruz, CA 95064, USA}




 \

\thanks{The work of Maria Schonbek was partially supported by NSF grant DMS-0600692}

\thanks{The work of Clayton Bjorland was partially supported by NSF grant OISE-0630623}

\begin{abstract}

This paper addresses the question of change of decay rate from exponential to algebraic for diffusive evolution equations. We show how the behaviour of the spectrum of the Dirichlet Laplacian in the two cases yields the passage from exponential decay in bounded domains to algebraic decay or no decay at all in the case of unbounded domains.  It is well known that such rates of decay exist: the purpose of this paper is to explain what makes the change in decay happen. We also discuss what kind of data is needed to obtain various decay rates.

\end{abstract}

\maketitle

\bigskip

\section{Introduction}

The purpose of this paper it to address the following two questions:

$\circ$  What makes solutions to diffusive evolution equations, with underlying linear part modelled by the heat equation, dramatically jump from exponential decay, when considered on a bounded domain $\Omega \subset    \mathbb{R}^n$, to algebraic decay or decay without a rate when considered on the whole of $\mathbb{R}^n$. This will be referred to the {\em decay-change phenomenon} (DCP).

$\circ$ What conditions are required on the data to ensure specific rates of energy decay?

\medskip

It is well known, that solutions to the heat equation (and solutions to similar linear second-order parabolic partial differential equations) defined on a bounded Lipschitz domain  $\Omega\subset \mathbb{R}^n$, $n \geq 1$, with initial datum $u_0$ in $\LL^2(\Omega)$, and subject to homogeneous Dirichlet boundary condition, decay exponentially in the $\LL^2(\Omega)$ norm. This is an easy consequence of the Poincar\'{e} inequality:

\[ \mbox{$\forall v \in \HH^1_0(\Omega)$}\qquad \lambda_1 \|v\|_2^2 \leq \|\nabla v\|_2^2, \]

where $\lambda_1>0$ is the smallest eigenvalue of the Dirichlet Laplacian on $\Omega$, defined by

\[-\Delta \,:\, v \in  \HH^1_0(\Omega) \mapsto -\Delta v \in \HH^{-1}(\Omega).\]

It turns out that  $\lambda_1 := C_0/d^2$, where $d= \mbox{diam}(\Omega)$ is the diameter of $\Omega$ and $C_0$ is a positive constant dependent only on the shape, but not the diameter of $\Omega$. This is easily seen by performing the change of variable $\Omega \subset \mathbb{R}^n \mapsto \hat{\Omega} = \frac{1}{d}(\Omega - x_0)\subset \mathbb{R}^n$ where $x_0$ is the barycenter of $\Omega$, and noting that $\mbox{diam}(\hat{\Omega})=1$. Thus, as  the domain grows, $\lim_{d \rightarrow \infty} \lambda_1 = 0$, and the Poincar\'{e} inequality is lost in the limit of $d \rightarrow \infty$.

Unsurprisingly, on $\mathbb{R}^n$ one can construct solutions to the homogeneous heat equation that decay at a very slow algebraic rate, and even ones that decay but without any rate. Moreover, we can find, for each time $T>0$  and each $\varepsilon \in (0,1)$, a solution $u$ to the homogeneous heat equation $u_t - \Delta u = 0$ on $\mathbb{R}^n$ with initial datum $u_0 \in V_{\beta} := \{ v: \| v \|_2^2 = \beta <\infty\}$ such that

\[ \frac{\|u(\cdot,T)\|}{\beta} \geq 1-\varepsilon.\]

One can also show the following theorem.

\begin{Theorem} \label{theo:norate}

There exists no function $G\,:\,(0,\infty)^2 \mapsto G(t,\beta)\in \mathbb{R}_+$ such that, if $u$ is a solution to the homogeneous heat equation $u_t - \Delta u = 0$ on $ \mathbb{R}^n$ with initial datum  $u_0 \in V_{\beta}$, then

\[\|u(\cdot,t)\|_2\leq G(t,\beta) \qquad \mbox{and}\qquad \lim_{t\rightarrow\infty}G(t,\beta)= 0 \quad \mbox{for all $\beta>0$}.
\]

\end{Theorem}

Analogous results hold for solutions to many nonlinear equations with a diffusive term modelled by the Dirichlet Laplacian or a fractional Dirichlet Laplacian, including the Navier--Stokes, Navier--Stokes-alpha, the quasi-geostrophic, and the magneto-hydrodynamics equations, (see, for example, \cite{CS}, \cite{MR837929}).

To explain the DCP, we shall take a careful look at the behavior of the spectrum of the Dirichlet Laplacian near zero on a bounded domain $\Omega\subset \mathbb{R}^n$. Although the distance to zero of the smallest eigenvalue of the Dirichlet Laplacian on $\Omega$ decays to zero as the diameter of $\Omega$ tends to infinity, for any fixed $\Omega$ the distance of the smallest eigenvalue to the origin remains positive, so the spectrum, which is entirely discrete and

consists of eigenvalues, remains bounded away from zero. This, as we will try to show, is the reason why there is no transition to slower than exponential decay, as $t \rightarrow \infty$, of $\|u(\cdot,t)\|_2$ for solutions $u$ of the homogeneous heat equation $u_t - \Delta u = 0$ in bounded domains $\Omega$, subject to homogeneous Dirichlet boundary condition.  We will clarify this issue by first considering the heat equation in a very simple one-dimensional situation (viz. on bounded intervals), and then extending the results to higher dimensions. For unbounded domains we will derive an extension of the Poincar\'{e} inequality, which will highlight the role played by neighborhoods of zero in frequency space.

This extension of the Poincar\'e inequality will be shown to be optimal in a sense that will be explained below.

The final sections will focus on finding the most general data for solutions to the heat equation in unbounded domains so that the corresponding solutions decay at a specific rate, then extend these ideas to certain nonlinear evolution equations whose solutions satisfy energy inequalities of the type

\[ \frac{{\rm d}}{{\rm d}t} \int _{\rn} |u(x,t)|^2 {\rm d}x \leq -C  \int _{\rn} |\nabla u(x,t)|^2 {\rm d}x. \]

\section{Notation}

In this section we collect the notation that will be used throughout the paper. The Fourier transform of $v\in \mathcal{S}(\rn)$ is defined by

\[\widehat v(\xi) = (2\pi)^{-n}\int_{\rn} {\rm e}^{-{\rm i}x\cdot\xi} v(x)\;{\rm d}x,\]

extended as usual to $\mathcal{S}'$. For a function $v:\rn\to \C$ and a
multi-index $\gm=(\gm_1,\gm_2,\dots,\gm_n)\in \mathbb{N}_0^n$, $D^\gm v$ denotes differentiation of
order $|\gm|=\gm_1+\gm_2 + \dots + \gm_n$ with respect to the $n$ (spatial) variables. If $v$ also
depends on time $t$, the symbol $D_t^jv$ is used to denote $j$th derivative of $v$ with respect to $t$.
If $k$ is a nonnegative integer, $\WW^{k,p}(\rn)$ will signify, as is standard, the Sobolev space consisting of functions in $\LL^p(\rn)$ whose generalized derivatives up to order $k$ belong to $\LL^p(\rn)$, $1\leq p \leq \infty$. When $p=2$, $\WW^{k,2}(\rn)=\HH^k(\rn)$, where the space $\HH^s(\rn)$ is defined for all $s\in \R$ as the space

of all $f\in \mathcal{S}'$ such that $(1+|\xi|^2)^{s/2}\hat f(\xi)\in \LL^2(\rn)$.

\section{Preliminaries}

In this section we recall some well-known facts regarding the basis of eigenfunctions for the Dirichlet Laplacian in bounded Lipschitz domains $\Omega \subset \mathbb{R}^n$, and their connection to the Poincar\'e inequality.

The next result is well known (see \cite{Evans}, for example).

\begin{Theorem}\label{theo:evans}

There is a countable orthonormal basis for $\LL^2(\Omega)$ which consists of eigenfunctions of the Dirichlet Laplacian $-\Delta\,:\, \HH^1_0(\Omega) \mapsto \HH^{-1}(\Omega)$.  The eigenfunctions belong to $\HH^1_0(\Omega)$
and the associated eigenvalues are all positive and bounded away from zero.

\end{Theorem}

To begin, we will suppose that $n=1$ and $\Omega= (-R,R)$.

\begin{Corollary}\label{sc}
The set of functions, where $\quad n=0,1,2,\dots $
\[\mathcal{B} = \left\{\sqrt{\frac{2}{R}}\sin\left(\frac{n\pi x}{R}\right),\qquad \sqrt{\frac{2}{R}}\cos\left(\frac{(2n+1)\pi x}{2R}\right),\right\}
\]
form an orthonormal basis for $\LL^2(-R,R)$.
\end{Corollary}

\begin{proof}

In the light of the previous theorem it suffices to seek $\lambda \in \mathbb{C}$ such that the two-point boundary-value problem

 \[\lambda w +\frac{{\rm d}^2 w}{{\rm d}x^2}  = 0, \quad x \in (-R,R),\qquad w(-R)=w(R)=0,\]

has a nontrivial solution $w \in \HH^1(-R,R)$. There is a countable set %

\[\{\lambda_n\,:\,n=1,2,\dots\}\cup\{\tilde{\lambda}_n\,:\,n=0,1,\dots\}\]

of $\lambda$'s for which such functions $w$ exist. An easy computation shows that

\[ w(x) = \left\{\begin{array}{ll} \sin\left(\frac{n\pi x}{R}\right) & \mbox{when $\lambda=\lambda_n := \left(\frac{n\pi }{R}\right)^2$},\quad n\in\{1,2,\dots\},\\

\cos\left(\frac{(2n+1)\pi x}{2R}\right) & \mbox{when $\lambda =\lambda_n := \left(\frac{(2n+1)\pi}{2R} \right)^2$},\quad n\in \{0,1,\dots\}.
\end{array}
\right.
\]

Hence, due to  Theorem \ref{theo:evans},  the collection of functions in the statement of the Corollary form a basis for $\LL^2(-R,R)$.

\end{proof}

\bigskip

Let $\{\lambda_n , \tilde{\lambda}_n\}_n$  be as defined above. As a consequence of Corollary \ref{sc},
any real-valued function $u \in \LL^2(-R,R)$
can be expanded in terms of a countable basis consisting of complex exponentials, as follows:

\bg \label{eq: sol}
u(x) = \frac{1}{2R}\left[\sum_{n=-\infty, n\neq 0}^{\infty}\widehat{u}(n) \exp \left({\rm i}\frac{n\pi}{R}x\right)   + \sum_{n=-\infty}^{\infty}
\widehat{\tilde{u}} (n) \exp\left( {\rm i} \frac{(2n+1)\pi}{2R} x\right)\right]
\ed

Here $\widehat{u}(n) = \widehat{u}(-n)$ and  $\widehat{\tilde{u}}(n) = -\widehat{\tilde{u}}(-n)$ (thus ensuring

that $u(x) \in \mathbb{R}$ for all $x \in (-R,R)$), where $\widehat{u}_n$, $n \in \mathbb{Z}\setminus \{0\}$,

and $\widehat{\tilde{u}}_n$, $n \in \mathbb{Z}$, are the corresponding  Fourier coefficients.

The main fact to note is that there is no term corresponding to $n=0$ in the first sum.

Let us now take an atomic measure $\nu$ supported on the set

\bg \label{M}
M= \left \{  \left\{ \frac{n\pi }{R}\right\}_n,\; n \in \mathbb{Z} \backslash \{0\};\quad\left\{\frac{(2n+1)\pi}{2R} \right\}_n,\; n\in \mathbb{Z} \right\},
\ed

and let $\mu_{\xi} = \frac{1}{2R}\nu$. Note that  $\min_{\xi \in M} |\xi| = \frac{\pi}{2R}$.
Then, the last expression (\ref{eq: sol}) can be rewritten as

\bg \label{eq:sol1}
u(x) = \int_{\mathbb{R}} \widehat{u}(\xi) \exp({\rm i} \xi x ) \; {\rm d}\mu_{\xi}  = \int_{M} \widehat{u}(\xi) \exp({\rm i} \xi x ) \; {\rm d}\mu_{\xi}.
\ed

We chose to write  $u$ as (\ref{eq:sol1}) so as to ensure
that we have a unified notation, regardless of whether we work on $\Omega = (-R,R)$ or the whole of $\mathbb{R}$.
Hence, by Parseval's identity,

\[ \|u\|_2^2 = \sum_{n=-\infty, n \neq 0}^{\infty}  |\widehat{u}(n) |^2 +  \sum_{n=-\infty}^{\infty}  |\widehat{\tilde{u}}(n) |^2
\]

Suppose now that $u \in \HH^1_0(-R,R)$.  Then,

\begin{align} 
 \frac{{\rm d}}{{\rm d}x} u (x) & =    \sum_{n=-\infty,  n\neq 0}^{\infty} \frac{{\rm i} n\pi}{R}\,\widehat{u}(n)\, \exp \left({\rm i}\frac{n\pi x}{R}\right) \nonumber \\ &  +  \sum_{n=-\infty }^{\infty} \frac{{\rm i} (2n+1)\pi }{2R}\,\widehat{\tilde{u}}(n)\, \exp\left({\rm i}\frac{(2n+1)\pi x}{R} \right)   \nonumber  .
 \end{align}

Hence the $\LL^2(-R,R)$ norm of the derivative is

  \begin{equation} \label{eq:der1}
  \left\|\frac{{\rm d}u}{{\rm d}x}\right\|_2^2 =
\sum_{n=-\infty, n\neq 0}^{\infty}  \left|\frac{n\pi}{R}\right|^2 \left|\widehat{u}(n)\right|^2\\ +
 \sum_{n=-\infty}^{\infty}  \left |\frac{(2n+1)\pi}{2R}\right |^2\left|\widehat{\tilde{u}}(n)\right|^2  \geq  \left(\frac{\pi}{2R}\right)^2\|u\|_2^2.
  \end{equation}

On rewriting the last inequality in (\ref{eq:der1}) in integral form, using the atomic measure $\mu_{\xi}$
defined above, we obtain

\bga\label{eq1}
\left\|\frac{{\rm d}u}{{\rm d}x}\right\|_2^2 = \int_{\{\xi:|\xi| \geq \frac{\pi}{2R}\} }|\xi|^2| \widehat{u}(\xi)|^2\; {\rm d} \mu_{\xi}
\geq  \left(\frac{\pi}{2R}\right)^2\|u\|_2^2.
\eda

\bigskip

\begin{Remark} The results we obtained above on $\Omega =(-R,R) \subset \mathbb{R}$ can be easily extended to
$\Omega =(-R,R)^n \subset \mathbb{R}^n$. The eigenfunctions are then products of sines and cosines. Hence we can, again, expand a real-valued function $u \in \LL^2(\Omega)$ into
a convergent (in $\LL^2(\Omega)$) infinite series of

complex exponentials. Upon doing so, we can express $u$ as the integral with respect to an atomic measure $\mu_\xi$ supported on the countable set that comprises the (discrete) spectrum of the Dirichlet Laplacian on $\Omega=(-R,R)^n$.

\bg \label{eq: sol1}
u(x) = \int_{\rn} \widehat{u}(\xi) \exp({\rm i} \xi \cdot x ) \; {\rm d}\mu_{\xi}  = \int_{M} \widehat{u}(\xi) \exp({\rm i} \xi \cdot x ) \; {\rm d} \mu_{\xi}.
\ed

As in the case of $n=1$, the support of the measure $\mu_\xi$ excludes $\xi=0$.\end{Remark}

\section{A Poincar\'{e}-like inequality on $\mathbb{R}^n$}

Poincar\'e's inequality is not valid when $\Omega = \mathbb{R}^n$. We shall nevertheless show that a modification of Poincar\'e's inequality holds on $\mathbb{R}^n$ by using Fourier transform instead of Fourier series

\begin{Theorem} \label{theo:UP} For each $u\in \HH^1(\rn)$ and any $\Lambda>0$, the following inequality holds:

\[\|\nabla u\|_2^2 \geq \Lambda^{2} \int_{\mathbb{R}^n} |\widehat{u}|^2 {\rm d}\xi  -  \int_{\{\xi:|\xi|\leq \Lambda\}} (\Lambda^2 - |\xi|^2) |\widehat{u}|^2 {\rm d}\xi.
\]

\end{Theorem}

\begin{proof}

This follows immediately by Plancherel's identity and by splitting the frequency domain as

\[\rn =  \mathcal{S} \cup \mathcal{S}^c \;\; \mbox{where}\;\; \mathcal{S} =\{ \xi :|\xi|\leq \Lambda\}, \]

that

\begin{eqnarray*}
\|\nabla u\|_2^2 &=& \int_{ \mathcal{S}^c} |\xi|^2 |\widehat{u}|^2 {\rm d}\xi +
\int_{ \mathcal{S}}|\xi|^2 |\widehat{u}|^2 {\rm d}\xi\\
&\geq&  \Lambda^2 \int_{ \mathcal{S}^c} |\widehat{u}|^2 {\rm d}\xi  +
\int_{ \mathcal{S}}|\xi|^2 |\widehat{u}|^2 {\rm d}\xi =  \Lambda^2 \int_{\rn} |\widehat{u}|^2 {\rm d}\xi - \int_{ \mathcal{S}} (\Lambda^2 - |\xi|^2) |\widehat{u}|^2 {\rm d}\xi,
\end{eqnarray*}

and the conclusion of the theorem follows.

\end{proof}

\bigskip

\begin{Remark}

Suppose that $\Omega$ is a bounded Lipschitz domain in $\mathbb{R}^n$ with Poincar\'e constant $C_\Omega$, that is,

\[ C_\Omega := \inf_{v \in \HH^1_0(\Omega)} \frac{\|\nabla v\|_2}{\|v\|_2}.\]

If in the last theorem we choose $\Lambda = \frac{\pi}{2R}$,  we recover the Poincar\'e inequality from the bounded domain $\Omega=(-R,R)$, minus an integral in frequency space over an interval of `radius' $\Lambda=  \left|\frac{\pi}{2R}\right|$ centred at $\xi=0$:

\[\|\nabla u(x)\|_2^2
\geq \left(\frac{\pi}{2R}\right)^2\|u\|_2^2 \; - \left(\frac{\pi}{2R}\right)^2\int_{\{\xi:|\xi|\leq \frac{\pi}{2R}\}}  |\widehat{u}|^2 {\rm d}\xi.\]

\end{Remark}

We now show that Poincar\'e's inequality for the whole space is optimal in the following sense:

\begin{Proposition}

Given any $K>0$ and any $\beta >0$, if there exist a $\mu>0$ such that

\[ \|u\|_2^2 \leq \frac{1}{K^2} \|\nabla u \|^2_2 +\mu\qquad \forall u \in \HH^1(\rn), \;\;\mbox{with}\;\; \|u(x)\|_2^2 =\beta\]

then

\[\mu \geq \int_{\{\xi :|\xi| \leq K \}} |\widehat{u}(\xi)|^2\; {\rm d}\xi= \Aa_{K}(u).\]

\end{Proposition}

\begin{proof}

We show by an example in $\R^2$ that given any $\mu$  we can find a function $u\in H^1(\R^2) \ni \|u(x)\|_2^2 =\beta $ for which the corresponding $ \Aa_{K}(u)$ is smaller then $\mu$. Similar examples can be found in any $R^n$ only that the computations are more tiresome.

We will work with the family of functions $u_{\alpha} =\frac{\sqrt{\beta}}{\sqrt{\pi}} \alpha e^{-[\frac{\alpha^2 |x|^2}{2}]}$.  Straightforward computations shows that

\bg \label{eq:1}
\|u_{\alpha} \|_2^2   = \beta \;\;\forall \alpha \;\;\mbox{and}\;\;\widehat{u}_{\alpha} =\frac{\sqrt{\beta}}{\sqrt{\pi}} \alpha^{-1}e^{-[\frac{\alpha^2 |\xi|^2}{2}]}
\ed

\bg \label{eq:2}
\|\nabla u_{\alpha} \|_2^2 = \beta \alpha^2  \;\;\mbox{and}\;\; \int_{\{\xi:|\xi|\leq K\}} | \widehat{u}_{\alpha}(\xi)|^2\;d\xi
= \beta(1-\e^{[\frac{K^2}{\alpha^2}]})
 \ed

Now using (\ref{eq:1}) and (\ref{eq:2}) in the Poincare inequality for the whole space gives

\bg \label{eq:3}
\beta \leq \frac{\alpha^2}{K^2} \pi + \beta(1-\e^{[\frac{K^2}{\alpha^2}]})
\ed

Let $ \beta (1-\e^{[\frac{K^2}{\alpha^2}]}) =\Bb(\alpha)$

We want to show that if

\bg \label{eq:4}
\beta \leq \frac{\alpha^2}{K^2} \beta + \mu
\ed

then $\mu \geq   \beta (1-\e^{[\frac{K^2}{\alpha^2}]})$

When $\alpha \geq K^2$  then the first term in the RHS in (\ref{eq:3}) is larger that the LHS of the inequality so $\mu=0$ will suffice to have the

\bg \label{eq:5}
\beta \leq \frac{\alpha^2}{K^2} \pi + \mu
\ed

We need to show that

\bg \label{eq:6}
\mu \geq  \Bb(\alpha) = \beta ( 1 - \frac{\alpha^2}{K^2}) \ed

for all $\alpha$,   in particular if we take  $\alpha = \frac{K^2}{2}$ it follows that it is necessary that $\mu \geq \frac{1}{2}$

Noting that

\[\lim_{\alpha \to 0} \Bb(\alpha) =0\]

it follows easily that we can find $\alpha$ close to zero so that  $\Bb(\alpha)\leq \frac{1}{2}$

\end{proof}

\subsection{Poincar\'e for Poisson type equations}

We now describe a different way we can write a Poincar\'e type inequality for the whole space. We will call this a ``Fake Poincar\'e inequalit'' (FPI)

\begin{Proposition} Let $u\in H^1(\rn)$ then there exist constants $K$ and $\alpha <1$ so that

\[ (1- \alpha) \|u\|_2^2 \leq K^{-2}  \|\nabla u\|_2^2\;\;\;\;\;(FPI) \]

\end{Proposition}

\begin{proof}

Use Poincar\'e for the whole space with $ \alpha = \frac{\int_{S_K}\widehat{u}\; d\xi}{\|u\|_2^2}$. Where
$S_k =\{\xi: |\xi |\leq K\} $ and we choose a K so that $\int_{S_K}\widehat{u}\; d\xi < \|u\|_2^2$.
\end{proof}

\bigskip

We now show how (FPI) can be used to obtain a Poincar\'e inequality for solutions to appropriate differential  equations.

\bigskip

\begin{Example}

Given some  constant $M$,

let  $\Vv_m^{\beta}=\{ f \in L^2:  |\widehat{f} | \leq M|\xi|^m,\;\; \mbox{for}\;\;|\xi|\leq  \beta\}$

Let $m\geq k$. Suppose $u$ satisfies

\begin{align}
D^{k}u = f\\
f \in Vv_m^{\beta},\;\;\mbox{some}\;\;m,\; \mbox{and}\; \beta \\
\lim_{|x| \to \infty} u =0\
\end{align}

then

 \[ \|u\|_2^2 \leq 2 \|\nabla u\|_2^2 \]

\end{Example}

\begin{proof}

Let $\Ss_{\alpha}= \{ \xi: |\xi| \leq \alpha\}$ and $M_o = M |\omega|$, where $|\omega|$ is the measure of the unit sphere.

Note that

\begin{itemize}
\item  $ \widehat{u} = \frac{\widehat{f}}{|\xi|^{2k}} $
\item  $\int_{\Ss_{\alpha}} |\widehat{u}|^2\; d\xi  = \int_{\Ss_{\alpha}} \frac{|\widehat{f}|^2}{|\xi|^{2k}}\; d\xi= M_o \beta^{2m-2k+n}$
\item Let $\beta_o $ be such that $ M_o \beta_o^{2m-2k+n} \leq \frac{\|f\|_2^2}{2}$
\end{itemize}

Then if we apply (FPI)   the conclusion of the example follows.

\end{proof}

\begin{Remark} Note that $\Vv_m^{\beta} \subset L^2 \cap I_m$, where $I_m =\{f: I_m(f) \in L^1(\rn)\}$,
with $I_m$ the Riesz potential of order $m$ for the function $f$.
\end{Remark}

\subsection{ Poincar\'e and decay}

Let  $u(x,t)\in H^1(\rn)$ satisfy

\bg \label{eq:FS} \frac{d}{dt} \|u(t)\|_2^2 \leq -C\|\nabla u\|_2^2\ed

The Poincar\'e  and the modified Poincar\'e inequalities  can be applied (for bounded or unbounded domains respectively) to the RHS of last inequality.

For  bounded domains to express the Fourier series as a Fourier integral, one uses a measure which has discrete support. Since  zero is not in the support the decay will the modified Poincar\'e inequality for the whole space that will induce the passage from exponential to algebraic decay or decay without a rate.

\medskip

To get decay on the whole space  for solutions $u$ of the inequality  (\ref{eq:FS}) we use Fourier splitting.

This method  applied  to the whole domain,  shows that we can look at  frequencies near the origin in balls that   depend on time. More precisely  Fourier splitting is just an application of the modified  Poincare  inequality in the whole space

where we have chosen $\Lambda =\Lambda(t)$  appropriately.

\begin{Theorem} \label{theo:NL} Let  $ u $ be a solution to

\begin{align} \label{eq:NL}
u_t = NL(u) + \Delta u\\
u_0\in L^2(\rn)
\end{align}

In addition suppose that the following properties hold

\begin{itemize}
\item $\int_{\rn} u \cdot NL(u) \;dx =0 $
\item $|\widehat{u}_o(\xi)| \leq C|\xi|^k \;\;\mbox{for}\;\;|\xi| << 1, \; \mbox{some}\;\;k\geq 0$
\item $NL(u)$ can be approximated in $L^2$ by  $NL(u_n)$ where $u_n$ are sufficiently smooth.
\end{itemize}

Then

 there exist a constant $C= C(u_o)$ so that

 \[\|u(t)\|_2 \leq C(t+1) ^{-\frac{n+2k}{4}}\]

 \end{Theorem}

\begin{proof}

Apply Poincar\'e for whole space and then Fourier Splitting to approximating sequences of solutions . See \cite{S} for this method. Then pass to the limit.

\end{proof}

\begin{Remark}

We note that the same theorem can be applied in case that diffusion is described by $(-\Delta)^{\alpha} u$ with $\alpha$ fractional, or by $D^m u$. The decay rates will have to be changes appropriately.
\end{Remark}

\section{ Classification of data}

As seen through Theorem \ref{theo:norate}, for arbitrary data just
in $L^2$ there is no specific decay rate for solutions of the heat
equation. That is for any fixed energy value we can find data which
leads to a solution whose heat energy decays arbitrarily fast or
slow.  Hence, the energy decay rate is dependent on the actual fom
of the data and not on initial energy.  In this section we
characterize the type of data which leads to exponential or
algebraic decay.

It is known that the structure of the data near the origin in
Fourier space dictates the rate of heat energy decay and our
theorems rely on this relation.  In particular we show that
exponential decay can occur if and only if the data is zero in some
neighborhood of the origin in Fourier space, we show that if a Riesz
potential of the initial data is in $L^1$ this can determine decay
rate, and finally we introduce a way to find what type of polynomial
``best'' describes the data near the origin and use this to
determine decay rates.  The last piece relies of finding a unique
decay character for any $L^2$ function by examining the Fourier
transform near the origin, this can be used to determine the
algebraic decay rate exactly.  We first will analyze solutions to
the Heat equation then extend the results to a more general setting
of parabolic equations which have a Laplacian linear part.

The starting idea is that functions with Fourier transform equal to zero near the origin decay exponentially, this suggests that a band
pass filter will be useful in characterizing such functions.

\begin{Lemma}

Given $\rho>0$, let $\chi_\rho(\xi)$ be the cut-off function equal to $1$ when $|\xi_j|\leq \rho_j, \;\;j=1,...,n$ and equal to zero elsewhere.  A function $u\in L^2(\mathbb{R}^n)$ satisfies $\hat{u}(\xi)=0$ for a.e. $|\xi_j|<\rho_j,\;\;j =1,...n$ if and only if $u = v - H_\rho\ast v$ a.e. for some $v\in L^2(\mathbb{R}^n)$ where $H_\rho(x)= \prod_{j=1}^nH_{\rho_j}(x_j)=\prod_{j=1}^n\frac{\sin{\rho_j x_j}}{x_j}$.  (Here $H$ is used to denote a \emph{high pass} filter).


\end{Lemma}

\begin{proof}

This proof is quickly checked by noting the Fourier transform of $H_{\rho}\ast u$ is $\chi_\rho\hat{u}$.

\end{proof}

The following theorem establishes that a solution to the Heat equation decays exponentially if and only if its initial data is zero a.e. in some possibly small ball centered at the origin in Fourier space.  Such functions have a particular form in the original space, shown by the previous lemma.  This characterizes all data which leads to exponential decay of heat energy.

\begin{Theorem}\label{thrm:expdecay}

The solution of the Heat equation with initial data $u_0\in L^2(\mathbb{R}^n)$ satisfies the decay bound  $\|u(t)\|_2^2\leq Ce^{-t\alpha^2}$ for some $\alpha>0$ and $C>0$ if and only if the initial data is of the form $u_0 =v_0-H_\rho\ast v_0$ for some $v_0\in L^2(\mathbb{R}^n)$.

\end{Theorem}

\begin{proof}

``$\Leftarrow$''

Let $\Bb =\{\xi: |\xi_j|<\rho_j\}$

Assume $u_0=v_0-H_\rho\ast v_0$, by the previous lemma $|\hat{u}(\xi)|^2=0$ for a.e. $\xi \in \Bb$.

\begin{align}
\|u(t)\|_2^2 &=\int_{\Bb}e^{-2|\xi|^2t}|\hat{u}_0|^2\, d\xi + \prod_{j=1}^n \int_{\Bb^c}e^{-2|\xi|^2t}|\hat{u}_0|^2\, d\xi\notag\\
&\leq \int_{\Bb}e^{-2|\xi|^2t}|\hat{u}_0|^2\, d\xi + e^{-2\rho^2t}\int_{\Bb^c}|\hat{u}_0|^2\, d\xi\notag\\
&\leq \int_{\Bb}e^{-2|\xi|^2t}|\hat{u}_0|^2\, d\xi + e^{-2\tilde{\rho}^2t}\|\hat{u}_0\|_2^2\notag
\end{align}

Where $\tilde{\rho} = \min\{\rho_j|j=1,...n\}$. By assumption the first integral on the RHS is zero.

\bigskip

``$\Rightarrow$''

By way of contradiction assume there is an $\alpha$ and a $C$ so that $\|u(t)\|_2^2 \leq Ce^{-t\alpha^2}$ and $\int_{|\xi|<\rho}|\hat{u}(\xi)|^2\, d\xi>0$ for all $\rho>0$.  Then, for $\rho=\frac{\alpha}{2}$, there is a $c>0$ so that $\int_{|\xi|<\rho}|\hat{u}(\xi)|^2\, d\xi>c$.  This implies:

\begin{align}
Ce^{-t\alpha^2}&\geq \|\hat{u}\|_2^2\notag\\
&= \int_{|\xi|<\rho}e^{-2|\xi|^2t}\hat{u}_0^2\, d\xi +\int_{|\xi|\geq\rho}e^{-2|\xi|^2t}\hat{u}_0^2\, d\xi\notag\\
&\geq e^{-2t\rho^2}\int_{|\xi|<\rho}\hat{u}_0^2\, d\xi\notag\\
&> ce^{-t\frac{\alpha^2}{2}}\notag
\end{align}

Taking $t$ sufficiently large violates this bound.

\end{proof}

The next Theorem attempts to characterize types of functions which lead to slower then exponential decay.  The first piece of this puzzle is the Riesz potential of initial data, $I_\beta(u_0)$.  It is defined in Fourier variables as

\begin{equation}
(\widehat{I_\beta u_0})(\xi)=\frac{\hat{u_0}(\xi)}{|\xi|^\beta}\notag
\end{equation}

Write, when the limit exists,

\begin{equation}\label{Apotential}
A_\beta(u_0) = \lim_{|\xi|\rightarrow 0}\frac{\hat{u}_0^2(\xi)}{|\xi|^\beta} = \int_{\mathbb{R}^n}(I_\beta u_0)(x)\,dx
\end{equation}

Note that this exists for all $I_{\beta} u_0\in L^1$.

\begin{Theorem}

Let $u$ be the solution of the heat equation corresponding to $u_0\in L^2$.


If  $I_{\beta} u_0\in L^1$ then

\begin{equation}
t^{\frac{n}{2}+\beta}\|\hat{u}(t)\|_2^2 \leq C (A_\beta(u_0))^2\notag
\end{equation}

Where  $A_\beta(u_0)$ was defined above. If $\|u(t)\|_2^2\leq C(1+t)^{\frac{n}{2}+\beta}$ for some $C$ and $\beta$ then

\begin{equation}
\liminf_{|\xi|\rightarrow 0}\frac{\hat{u}_0(\xi)}{|\xi|^\beta}<\infty\notag
\end{equation}
\end{Theorem}




\begin{proof}

Both statements in the theorem are consequences of the following equality, proved by the change of variables $\sqrt{t}\xi=\eta$:

\begin{align}
t^{\frac{n}{2}+\beta}\|\hat{u}(t)\|_2^2 &= t^{\frac{n}{2}+\beta}\int_{\mathbb{R}^n}e^{-2|\xi|^2t}\hat{u}_0^2(\xi)\,d\xi\notag\\
&=\int_{\mathbb{R}^n}e^{-2\eta^2}\eta^{2\beta}(\widehat{I_\beta w}(\eta))^2\, d\eta\notag
\end{align}

If $\lim_{|\xi|\rightarrow 0}\frac{\hat{u}_0(\xi)}{|\xi|^\beta}$ exists and is finite then the Lebesque dominated convergence theorem proves the first statement.  If $\|u(t)\|_2^2\leq C(1+t)^{-\frac{n}{2}-\beta}$ then

\begin{equation}
\int_{\mathbb{R}^n}e^{-2\eta^2}\eta^{2\beta}(\widehat{I_\beta w}(\eta))^2\, d\eta\notag \leq C\frac{t^{\frac{n}{2}+\beta}}{(1+t)^{\frac{n}{2}+\beta}}\notag
\end{equation}

Fatou's lemma then proves the second statement.

\end{proof}

The above theorem showns when $A_\beta(u_0)$ exists we can determine the
rate of decay from $\beta$ and the dimension of space but the
condition $A_\beta(u_0)$ exists is stronger then $u_0\in L^2$ so a
more general structure is needed.  Moreover, determining if
$A_\beta(u_0)$ exists entails determining if $I_{\beta}u_0 \in L^1$,
which might not be always simple.  On the other hand, if the Fourier transform of
initial data is a polynomial of the form $|\xi|^q$ in some possibly
small neighborhood of the origin then the solution will decay as
$~C(1+t)^{-q-\frac{n}{2}}$.  This can be checked by calculating (or
estimating) the integral $\int_{B(\rho)}
e^{-2|\xi|^2t}|\xi|^{2q}d\xi$.  Unfortunately only a small amount of
initial data can be described in this way and we are pushed to find
what order of polynomial \emph{best} describes a general function
$u_0\in L^2$ near the origin in Fourier space and using this
information to determine the decay rate.

\begin{Definition}

The ``decay indicator'' $P_q(u_0)$ is defined as follows.  Let $B(\rho)$ be the ball of radius $\rho$ centered at the origin, for $q\in(-\frac{n}{2},\infty)$:

\begin{equation}
P_q(u_0) = \lim_{\rho\rightarrow 0}\rho^{-2q-n}\int_{B(\rho)}|\hat{u}_0(\xi)|^2\, d\xi\notag
\end{equation}

\end{Definition}

When $u_0\in L^2(\mathbb{R}^n)$ the integral $\int_{B(\rho)}|\hat{u}_0(\xi)|^2\, d\xi$ considered as a function of $\rho$ is continuous, monotone decreasing, and bounded from below as $\rho$ becomes small, this is enough to ensure that $P_q$ is defined for all $q$ and $u_0\in L^2$.  It is trivially zero for all $q\leq-\frac{n}{2}$ when $u_0\in L^2$ so we consider only $q\in (-\frac{n}{2},\infty)$.  $P_q(u_0)$ compares $\hat{u}_0$ to the polynomial $|\xi|^{2q}$ near the origin and takes values in the (non-negative) extended real numbers, we are interested in three outcomes: $P_q(u_0) = 0,\infty,c$ with $c\neq 0$.  Depending on the outcome we think of $u_0$ as a polynomial with degree, respectively, greater, less, or equal to $q$ near the origin.  It is easy to check that $P_q(|\xi|^q)=|\omega|(n+2q)=\mu_q$.

Recalling (\ref{Apotential}), when $A_q(u_0)$ exists it bounds $P_q(u_0)$, in this sense $P_q$ is weaker then $A_q$.  Indeed,

\begin{align}
P_q(u_0)&= \lim_{\rho\rightarrow 0}\rho^{-2q-n}\int_{B(\rho)}|\hat{u}_0|^2\,d\xi\notag\\
&\leq \lim_{\rho\rightarrow 0} \left(\sup_{B(\rho)}\frac{|u_0|^2}{|\xi|^{2q}}\right)\rho^{-2q-n}\int_{B(\rho)}|\xi|^{2q}\,d\xi\notag\\
&=\lim_{\rho\rightarrow 0} \sup_{B(\rho)}\frac{|u_0|^2}{|\xi|^{2q}}\mu_q\notag\\
&=A^2_{q}((u_0))\mu_q\notag
\end{align}

\begin{Definition}

For a given $u_0\in L^2(\mathbb{R}^n)$, we call the unique value $q^*$ given by Lemma \ref{lemma:qstar} below  the ``decay character.''

\end{Definition}

\begin{Lemma}\label{lemma:qstar}

For each $u_0\in L^2$ there is at most one value of $q\in (-\frac{n}{2},\infty)$ so that $0<P_q(u_0)<\infty$.  If such a number exists we denote it $q^*$.  If no such number exits we take $q^*=-\frac{n}{2}$ in the case where $P_q(u_0)=0$ for all $q\in (-\frac{n}{2},\infty)$ and $q^*=\infty$ if $P_q(u_0)=\infty$ for all $q\in (-\frac{n}{2},\infty)$.

\end{Lemma}

\begin{proof}

Let $a = \sup\{q: P_q(u_0)=0\}$ and $b = \inf \{q: P_q(u_0)=\infty\}$.  If $q$ is such that $P_q(u_0)=\infty$ and $r>q$ then $P_r(u_0)=\infty$, this is observed by taking the limit of the following inequality which is valid for all $\rho<1$:

\begin{align}
\rho^{-2r-n}\int_{B(\rho)}|\hat{u}_0(\xi)|^2\, d\xi &= \rho^{2q-2r}\rho^{-2q-n}\int_{B(\rho)}|\hat{u}_0(\xi)|^2\, d\xi\\
&>\rho^{-2q-n}\int_{B(\rho)}|\hat{u}_0(\xi)|^2\, d\xi\notag
\end{align}

From this we conclude $P_q(u_0)=\infty$ for all $q>b$ and the lemma is true if $b=-\frac{n}{2}$, a similar statement can be made concerning $q<a$ and shows the lemma is true if $a=\infty$.

It is also clear that $a\leq b$ and $q^*$ will be well defined for all $u_0\in L^2$ when the proof is finished.  To accomplish this we wish to show $a=b$, this is proved by contradiction.  Assume, contrary to the statement, there exists $q\in(a,b)$, then $P_q(u_0)=c$ for some $0<c<\infty$.  If $\epsilon>0$, similar to the above inequality:

\begin{align}
P_{q+\epsilon}(u_0) &= \left(\lim_{\rho\rightarrow 0}\rho^{-2\epsilon}\right)c\notag
\end{align}

This shows $P_{q+\epsilon}(u_0)=\infty$, since $\epsilon$ was chosen arbitrarily we conclude $q=b$.

\end{proof}

The decay character is calculated from the behavior of an $L^2$ function near the origin in Fourier space, we now prove a theorem relating decay rates and $P_q(u_0)$.  A consequence of this theorem (Theorem \ref{cor:decaycharacter}) will summarize the relation between the decay character and the exact decay rate which is our goal.

\begin{Theorem}\label{theorem:decaybounds}

Let $u$ be the solution to the heat equation corresponding to $u_0\in L^2(\mathbb{R}^n)$ and $q\in (-\frac{n}{2},\infty)$.  If $P_q(u_0)>0$ there exists a constant $C_1>0$ which depends only on $\|u_0\|^2_2$ and the dimension of space so that

\begin{equation}
C_1(1+t)^{-q-\frac{n}{2}} \leq \|u(t)\|_2^2 \notag
\end{equation}

If $P_q(u_0)<\infty$ there are constants $C_2,C_3>0$, again depending only on $\|u_0\|_2^2$ the dimension of space so that

\begin{equation}
\|u(t)\|_2^2 \leq C_2(C_3+t)^{-q-\frac{n}{2}}\notag
\end{equation}

\end{Theorem}

\begin{proof}

We consider first the lower bound assuming $P_q(u_0)>0$.   Relying on the existence of the limit and the fact that it is bounded away from zero we may take $\rho_0>0$ sufficiently small to ensure, for all $\rho\leq\rho_0$:

\begin{equation}
c_1<\rho^{-2q-n}\int_{B(\rho)}|\hat{u}_0|^2\,d\xi\notag
\end{equation}

Let $0<\rho(t)\leq\rho_0$, we will soon chose it exactly.

\begin{align}
\|u(t)\|_2^2 &= \int_{B(\rho)}e^{-2|\xi|^2t}|\hat{u}_0|^2\,d\xi +\int_{B^C(\rho)}e^{-2|\xi|^2t}|\hat{u}_0|^2\,d\xi\notag\\
&\geq (e^{-2\rho^2t}\rho^{-2q-n})\left(\rho^{-2q-n}\int_{B(\rho)}|\hat{u}_0|^2\,d\xi\right) \notag\\
&\geq (e^{-2\rho^2t}\rho^{2q+n})c_1\notag
\end{align}

Choosing $\rho(t)= \rho_0(1+t)^{-\frac{1}{2}}$ proves the lower bound.

To prove the upper bound assume $P_q(u_0)<\infty$ and take $\rho_0>0$ sufficiently small so that for all $\rho\leq\rho_0$:

\begin{equation}
\rho^{-2q-n}\int_{B(\rho)}|\hat{u}_0|^2\,d\xi\leq c_2\notag
\end{equation}

The constant $c_2$ is known to exist since $P_q(u_0)<\infty$.  We now use the Fourier Splitting Method (\cite{S}), starting with the well known energy inequality for the heat equation with $0<\rho(t)\leq\rho_0$.

\begin{align}
\frac{1}{2}\frac{d}{dt}\|u(t)\|_2^2&\leq -\|\nabla u(t)\|_2^2\notag\\
&\leq -\rho^2\int_{B^C(\rho)}|\hat{u}(t)|^2\,d\xi\notag
\end{align}

This implies

\begin{align}
\frac{d}{dt}\|u(t)\|_2^2 +2\rho^2\|u(t)\|_2^2 &\leq 2\rho^2\int_{B(\rho)}|\hat{u}(t)|^2\,d\xi\notag\\
&\leq 2\rho^{2+2q+n}C_2\notag
\end{align}

Take $m>q+\frac{n}{2}$ and choose $\rho(t)=\frac{m}{2(C_3+t)^{1/2}}$ where $C_3>0$ is large enough to guarantee $\rho(0)\leq \rho_0$.  Solve the differential inequality with the integrating factor $(C_3+t)^m$ to find

\begin{equation}
\|u(t)\|_2^2\leq C(C_3+t)^{-q-\frac{n}{2}}+(C_3+t)^{-m}\|u_0\|_2^2\notag
\end{equation}

This is the upper bound.

\end{proof}

\begin{Theorem}\label{cor:decaycharacter}

Let $u_0\in L^2(\mathbb{R}^n)$, $u(t)$ the corresponding solution to the heat equation, and $q^*$ the decay character.  If $-\frac{n}{2}<q^*<\infty$ then there are constants $C_1,C_2,C_3>0$ so that

\begin{equation}
C_1(1+t)^{-q^*-\frac{n}{2}} \leq \|u(t)\|_2^2 \leq C_2(C_3+t)^{-q^*-\frac{n}{2}}\label{bounds1}
\end{equation}

Moreover, if $q^*=-\frac{n}{2}$ or $q^*=\infty$ then $\|u(t)\|_2^2$ decays, respectively, slower or faster then  $(1+t)^{-q-\frac{n}{2}}$.

\end{Theorem}

\begin{proof}

If $-\frac{n}{2}<q^*<\infty$ then $0<P_{q^*}(u_0)<\infty$ and (\ref{bounds1}) follows from the previous theorem.  If $q^*=-\frac{n}{2}$ we have, for all $q\in(-\frac{n}{2},\infty)$, by the previous theorem we are guaranteed the existence of a constant $C_1(q)>0$ such that

\begin{equation}
C_1(1+t)^{-q-\frac{n}{2}} \leq \|u(t)\|_2^2\notag
\end{equation}

Letting $q$ take all values in $(-\frac{n}{2},\infty)$ proves the slow decay statement.  The statement concerning $q^*=\infty$ is argued similarly.

\end{proof}

\section{Decay for solutions to some nonlinear evolution equations}

Our next goal is to apply this idea to a more general class of PDEs
with a non-linear term.  The idea is that if
the non-linear term decays faster then the linear term, the decay of
solutions will be given through the decay character.  If the
non-linear term decays slower then it will determine an upper bound
on the decay rate, a lower bound will require more knowledge of the
specific non-linear structure.  First we will assume a bound
on the non-linear term in Fourier space and derive a bound on
energy decay.  After this we will use results from \cite{W} with the
decay character to describe energy decay of solutions for the
Navier-Stokes equation.

\begin{Theorem}

Let $u$ be a solution of the PDE (\ref{eq:NL})and assume the non-linearity satisfies the following conditions:

\begin{itemize}

\item[1.] We are justified in writing the solution as

\begin{equation}\notag
\hat{u}= e^{-|\xi|^2t}\hat{u_0} + \int_0^te^{-|\xi|^2(t-s)}\widehat{NL(u)}(\xi,s)\,ds
\end{equation}

\item[2.] $\widehat{NL(u)}(\xi,s) \leq C|\xi|^k$

\item[3.] $\int_{\rn} u \cdot NL(u) \;dx =0$

\end{itemize}

Let $q^*$ be the decay character associated with $u_0$, then for any $q<q^*$ there exist constants $C_1,C_2>0$ so the solution of the PDE satisfies the energy decay estimate

\begin{equation}
\|u(t)\|_2^2\leq C_1(C_2+t)^{-q-\frac{n}{2}}+C_1(C_2+t)^{-k-\frac{n}{2}}\notag
\end{equation}

\end{Theorem}

\begin{Remark} The assumptions in this theorem are reasonable for equations such as
the Navier-Stokes equation, the Navier-Stokes-$\alpha$ equation and
the Magneto-Hydrodynamic equation among others (see, for example,
\cite{CS}, \cite{S}, and \cite{MR1380452}). The decay rate in this
theorem is determined by $q^*$ or $k$, whichever is smaller.
\end{Remark}

\begin{proof}

Assumptions 1 and 2 imply

\begin{equation}
|\hat{u}|\leq |e^{-|\xi|^2t}\hat{u_0}| + C\int_0^te^{-|\xi|^2(t-s)}|\xi|^k\,ds\notag
\end{equation}

Completing the integral on the RHS then squaring yields

\begin{equation}
|\hat{u}|^2\leq e^{-2|\xi|^2t}|\hat{u_0}|^2 + C|\xi|^{2k-2}\notag
\end{equation}

Consider now any $q\leq q^*$ and take $\rho_0>0$ sufficiently small so that for all $\rho\leq\rho_0$

\begin{equation}
\rho^{-2q-n}\int_{B(\rho)}|\hat{u}_0|^2\,d\xi\leq C_3\notag
\end{equation}

Here $C_3$ is some constant known to exist since $q<q^*$.  Assumption 3 allows an energy inequality from which to use the Fourier Splitting Method, proceeding now as in the proof of Theorem \ref{theorem:decaybounds}:

\begin{align}
\frac{d}{dt}\|u(t)\|_2^2 +2\rho^2\|u(t)\|_2^2 &\leq 2\rho^2\int_{B(\rho)}|\hat{u}(t)|^2\,d\xi\notag\\
&\leq \rho^2\int_{B(\rho)}|\hat{u}_0|^2\,d\xi + C\rho^2\int_{B(\rho)}|\xi|^{k-2}\,d\xi\notag\\
&\leq 2\rho^{2+2q+n}C_3 + C\rho^{2k+n}\notag
\end{align}

Take $m>\max\{q+\frac{n}{2},k+\frac{n}{2}\}$ and choose $\rho(t)=\frac{m}{2(C_2+t)^{1/2}}$ where $C_2>0$ is large enough to guarantee $\rho(0)\leq \rho_0$. To finish the Fourier Splitting Method multiply the above equation by $(C_2+t)^m$ and solve the differential inequality to find

\begin{equation}
\|u(t)\|_2^2\leq C(C_2+t)^{-q-\frac{n}{2}}+C(C_2+t)^{-k-\frac{n}{2}}+(C_2+t)^{-m}\|u_0\|_2^2\notag
\end{equation}

This finishes the proof.

\end{proof}

In the specific case of the Navier-Stokes equation there has been
significant investigation into the rate at which a solution
approaches a solution of the heat equation with the same initial
data, see \cite{Carpio}, \cite{Miyakawa}, \cite{W}.  We will now
demonstrate how the results of this section fit with the main theorem in
\cite{W}.

\begin{Theorem}(Wiegner)\label{theorem:wiegner}

Let $n\geq 2$ be a weak solution of the Navier-Stokes equation on $\mathbb{R}^n$, $2\leq n\leq 4$,

\begin{align}\label{PDE:NS}
\partial_t u-\triangle u+u\cdot\nabla u +\nabla p=0\\
\nabla\cdot u=0\ \ \ \ u(0)=u_0\notag
\end{align}

which satisfies the energy inequality

\begin{equation}
\|u(t)\|_2^2+2\int_s^t\|\nabla u(r)\|_2^2\,dr\leq \|u(s)\|^2_2\notag
\end{equation}

If $\|e^{\triangle t}u_0\|_2^2\leq C(1+t)^{-\alpha}$, then $\|u(t)-e^{\triangle t}u_0\|_2^2\leq h_{\alpha}(t)(1+t)^{-d}$ with $d=\frac{n}{2}+1-2\max\{1-\alpha,0\}$  and

\bg
 h_{\alpha}(t)=\left\{
 \begin{array}{lr}
 \epsilon(t) &\;\;\mbox{ for} \ ;\;\alpha=0,\mbox{with} \;\;\epsilon(t)\searrow 0\ for\ t\rightarrow \infty\notag\\
C ln^2(t+e)&\;\; \mbox{for} \;\; \alpha=1\notag\\
C&\;\;\mbox{ for} \;\;\; \alpha\neq 0,1
\end{array} \right.
\ed






\end{Theorem}

\begin{Remark}
The actual theorem proved by Wiegner is more general as it includes a forcing term, forcing terms are outside the scope of this paper.
\end{Remark}

\begin{proof}

See \cite{W}

\end{proof}

\begin{Theorem}\label{cor:NSdecaycharacter}
Let $u$ be as in Theorem \ref{theorem:wiegner} and $q^*$ be the decay character associated with $u_0$.

\begin{itemize}
\item[1.]If $-\frac{n}{2}<q^*< 1-\frac{n}{2}$ then there are constants $C_1,C_2,C_3>0$ so that
\begin{equation}
C_1(1+t)^{-q^*-\frac{n}{2}} \leq \|u(t)\|_2^2 \leq C_2(C_3+t)^{-q^*-\frac{n}{2}}\notag
\end{equation}
\item[2.]If $q^*\geq 1-\frac{n}{2}$ then $\|u(t)\|_2^2\leq C(1+t)^{-\frac{n}{2}-\beta}$ where $\beta=\min(q^*,1)$.
\item[3.]If $q^*=-\frac{n}{2}$ and $n=3,4$ then $\|u(t)\|_2^2$ decays slower then any polynomial.
\end{itemize}

\end{Theorem}

\begin{proof}

This is a combination of the the above theorem and Theorem \ref{cor:decaycharacter}.  In case 1, Corollary \ref{cor:decaycharacter} gives

\begin{equation}
C_1(1+t)^{-\frac{n}{2}-q^*}\leq \|e^{\triangle t}u_0\|_2^2\leq C_2(C_3+t)^{-\frac{n}{2}-q^*}\notag
\end{equation}

while Theorem \ref{theorem:wiegner} allows ($\alpha = q^*+\frac{n}{2}$, $d=\frac{3n}{2}-1-2q^*$):

\begin{equation}
\|u(t)-e^{\triangle t}u_0\|_2^2\leq C(1+t)^{-\frac{3n}{2}+1-2q^*}\notag
\end{equation}

Notice that the second decays faster.  Combining these with the triangle inequality proves 1.

Case 2 is similar, this time Theorem \ref{cor:decaycharacter} gives

\begin{equation}
C_1(1+t)^{-\frac{n}{2}-q^*}\leq \|e^{\triangle t}u_0\|_2^2\leq C_2(C_3+t)^{-\frac{n}{2}-q^*}\notag
\end{equation}

while Theorem \ref{theorem:wiegner} allows ($\alpha=\frac{n}{2}+q^*$, $d=\frac{n}{2}+1$):

\begin{equation}
\|u(t)-e^{\triangle t}u_0\|_2^2\leq h_\alpha(t)(1+t)^{-\frac{n}{2}-1}\notag
\end{equation}

An application of the triangle inequality shows

\begin{equation}\notag
\|u(t)\|_2^2 \leq h_\alpha(t)(1+t)^{-\frac{n}{2}-1}+C_2(C_3+t)^{-\frac{n}{2}-q^*}
\end{equation}

Case 3 also follows from Theorems \ref{cor:decaycharacter} and

\ref{theorem:wiegner}, and the triangle inequality.  With $\alpha
=0$ in Theorem \ref{theorem:wiegner} we have $\|u(t)-e^{\triangle
t}u_0\|_2^2\leq \epsilon(t)(1+t)^{-\frac{n}{2}+1}$ while
$\|e^{\triangle t}u_0\|_2^2$ decays slower then any polynomial
(recall $\epsilon(t)\searrow 0$). When $n=3,4$ the first decays
faster, thus, for large $t$,

\begin{equation}\notag
2\|e^{\triangle t}u_0\|_2^2 \leq |\|e^{\triangle t}u_0\|_2^2-\|u(t)-e^{\triangle t}u_0\|^2_2| \leq \|u(t)\|_2^2
\end{equation}

\end{proof}

\bibliographystyle{plain}

\bibliography{heat-stuff-cm}

\begin{thebibliography}{1}

\bibitem{CS}
C.~Bjorland and M.~E. Schonbek.
\newblock On questions of decay and existence for the viscous camassa-holm
  equations.

\bibitem{Carpio}
A.~Carpio.
\newblock Asymptotic behavior for the vorticity equations in dimensions two and
  three.
\newblock {\em Comm. Partial Differential Equations}, 19(5-6):827--872, 1994.

\bibitem{Evans}
L~.~C. Evans.
\newblock Partial differential equations.

\bibitem{Miyakawa}
T.~Miyakawa and M.~E. Schonbek.
\newblock On optimal decay rates for weak solutions to the navier-stokes
  equations in $\mathbb{R}^n$.
\newblock {\em Mathematica Bohemica}, 126(2):443--455, 2001.

\bibitem{S}
M.~E. Schonbek.
\newblock {$L\sp 2$} decay for weak solutions of the {N}avier-{S}tokes
  equations.
\newblock {\em Arch. Rational Mech. Anal.}, 88(3):209--222, 1985.

\bibitem{MR837929}
M.~E. Schonbek.
\newblock Large time behaviour of solutions to the {N}avier-{S}tokes equations.
\newblock {\em Comm. Partial Differential Equations}, 11(7):733--763, 1986.

\bibitem{MR1380452}
M.~E. Schonbek, T.~P. Schonbek, and Endre S{\"u}li.
\newblock Large-time behaviour of solutions to the magnetohydrodynamics
  equations.
\newblock {\em Math. Ann.}, 304(4):717--756, 1996.

\bibitem{W}
M.~Wiegner.
\newblock Decay results for weak solutions of the {N}avier-{S}tokes equations
  on {${\bf R}\sp n$}.
\newblock {\em J. London Math. Soc. (2)}, 35(2):303--313, 1987.

\end{thebibliography}

\end{document}